\documentclass[12pt]{article}
\usepackage{amsfonts}

\usepackage{mathrsfs}
\usepackage{amscd}
\usepackage{amsmath,amsfonts,amssymb,amscd}
\usepackage{indentfirst,graphics,epsfig,psfrag}
\input{epsf}

\usepackage{amsmath,amssymb,amsthm, amscd, amsfonts, graphicx,ifpdf}
\usepackage{lineno}

\newtheorem{thm}{Theorem}[section]

\newtheorem{lem}{Lemma}[section]
\newtheorem{conj}[thm]{Conjecture}

\newtheorem{obv}[thm]{Observation}
\def\qed{\nopagebreak\hfill{\rule{4pt}{7pt}}
\medbreak}

\def\pf{\noindent {\it Proof.} }

\title{\bf\Large $k$-tree connectivity of line graphs}
\author{
\small Shasha Li\\
\small School of Mathematics and Statistics,
Ningbo University, Ningbo 315211, Zhejiang, China\\
\small Email:  lishasha@nbu.edu.cn\\
}
\date{}
\begin{document}

\maketitle

\begin{abstract}
For a graph $G=(V,E)$ and a set $S\subseteq V(G)$ of size at least $2$,
an {\it $S$-Steiner tree} $T$ is a subgraph of $G$ that is a tree with
$S\subseteq V(T)$. Two $S$-Steiner trees $T$ and $T'$ are {\it internally
disjoint} (resp. {\it edge-disjoint}) if $E(T)\cap E(T')=\emptyset$ and $V(T)\cap V(T')=S$
(resp. if $E(T)\cap E(T')=\emptyset$). Let $\kappa_G (S)$ (resp. $\lambda_G (S)$) denote the maximum number
of internally disjoint (resp. edge-disjoint) $S$-Steiner trees in $G$. The
{\it $k$-tree connectivity} $\kappa_k(G)$ (resp. {\it $k$-tree edge-connectivity} $\lambda_k(G)$) of $G$ is then defined
as the minimum $\kappa_G (S)$ (resp. $\lambda_G (S)$), where $S$ ranges over all $k$-subsets of $V(G)$.
In [H. Li, B. Wu, J. Meng, Y. Ma, Steiner tree packing number and tree
connectivity, Discrete Math. 341(2018), 1945--1951], the authors
conjectured that if a connected graph $G$ has at least $k$ vertices and at least $k$ edges,
then $\kappa_k(L(G))\geq \lambda_k(G)$ for any $k\geq 2$, where $L(G)$ is the line graph of $G$. In this paper,
we confirm this conjecture and prove that the bound is sharp.
\end{abstract}

\section{Introduction}
The connectivity $\kappa(G)$ of a graph $G$ is the minimum cardinality of
a subset $V'$ of vertices such that $G-V'$ is disconnected or trivial.
The edge-connectivity $\lambda(G)$ of a graph $G$ is the minimum cardinality of
a subset $E'$ of edges such that $G-E'$ is disconnected.
An equivalent definition of connectivity was given in \cite{Whitney}.
For each $2$-subset $S=\{u,v\}$ of vertices of $G$,
let $\kappa(S)$ denote the maximum number of internally vertex-disjoint
$(u,v)$-paths in $G$. Then $\kappa(G)=$min$\{\kappa (S)|S\subseteq V \text{ and } |S|=2\}$.
Similarly, the edge-connectivity also has two equivalent definitions.
Let $\lambda(S)$ denote the maximum number of edge-disjoint
$(u,v)$-paths in $G$. Then $\lambda(G)=$min$\{\lambda (S)|S\subseteq V \text{ and } |S|=2\}$.

As a means of strengthening the connectivity, the tree connectivity was introduced
by Hager \cite{Hager1,Hager2} (or generalized connectivity
by Chartrand et al. \cite{Chartrand2}) to meet wider applications.
Given a graph $G=(V,E)$ and a set $S\subseteq V(G)$ of size at least $2$,
an {\it $S$-Steiner tree} or a {\it Steiner tree connecting $S$}
is such a subgraph $T=(V', E')$ of $G$ that is a tree with
$S\subseteq V'$. Two $S$-Steiner trees $T$ and $T'$ are said to be {\it internally
disjoint} if $E(T)\cap E(T')=\emptyset$ and $V(T)\cap V(T')=S$.
Let $\kappa_G (S)$ denote the maximum number
of internally disjoint $S$-Steiner trees in $G$. The {\it $k$-tree connectivity}
(or {\it generalized $k$-connectivity}) of $G$, denoted by $\kappa_k(G)$, is then defined
as $\kappa_k(G)=$min$\{\kappa_G(S)| S\subseteq V(G)\ and\ |S|=k\}$,
where $2\leq k\leq n$. Clearly, when $k=2$, $\kappa_2(G)$ is
exactly the classical connectivity $\kappa(G)$.

As a natural counterpart of the tree-connectivity, the tree edge-connectivity
(or generalized edge-connectivity) was introduced by Li et al. \cite{LMS}.
For $S\subseteq V(G)$ and $|S|\geq 2$, let $\lambda_G (S)$ denote the maximum number
of edge-disjoint $S$-Steiner trees in $G$. The {\it $k$-tree edge-connectivity}
(or {\it generalized $k$-edge-connectivity}) of $G$, denoted by $\lambda_k(G)$, is then defined
as $\lambda_k(G)=$min$\{\lambda_G(S)| S\subseteq V(G)\ and\ |S|=k\}$,
where $2\leq k\leq n$. It is also clear that when $k=2$, $\lambda_2(G)=\lambda (G)$.

There have been many results on the
$k$-tree (edge-)connectivity, see \cite{Chartrand,LLZ,LIM,LMS,LinZ}
and a book \cite{LIMAO}.

The {\it line graph} $L(G)$ of $G$ is the graph whose
vertex set can be put in one-to-one correspondence with the edge set of $G$ in
such a way that two vertices of $L(G)$ are adjacent if and only if the corresponding
edges of $G$ are adjacent. The connectivity of the line graph of a graph $G$
is closely related to the edge-connectivity of $G$.

\begin{lem}[Chartrand and Stewart \cite{Chartrand3}]\label{lem1}
If $G$ is a connected graph, then $\kappa(L(G))\geq \lambda(G)$.
\end{lem}

Naturally, one would like to study the relationship between $\kappa_k(L(G))$ and $\lambda_k(G)$,
for $k\geq 3$. In \cite{LMS}, Li et al. showed that
if $G$ is a connected graph, then $\kappa_3(L(G))\geq \lambda_3(G)$.
In \cite{LWMM}, Li et al. showed that
if a graph $G$ is connected, then $\kappa_4(L(G))\geq \lambda_4(G)$.
Furthermore, they proved that if a connected graph $G$ has at least $k$
vertices and at least $k$ edges, then $\kappa_k(L(G))\geq \lfloor\frac{3\lambda_k(G)}{4}\rfloor-1$
for any $k\geq 2$. However, they suspect that their result is not sharp
and proposed the following conjecture:

\begin{conj}[Li, Wu, Meng and Ma \cite{LWMM}]\label{conj1}
Let $k\geq 2$ be an integer. If a connected graph $G$ has at least $k$ vertices and at least $k$ edges, then
$\kappa_k(L(G))\geq \lambda_k(G)$.
\end{conj}

In this paper, we will confirm this conjecture and prove that the bound is sharp.

\section{Main result}

Before proving our main result, we first introduce some concepts.
A maximal connected subgraph of $G$ is called a {\it component} of $G$.
A connected acyclic graph is called a $tree$.
The vertices of degree 1 in a tree are called $leaves$.
A connected graph $G$ with $|V(G)|=|E(G)|$ is called a {\it unicyclic graph}.
A $spanning$ $subgraph$ of a graph $G$ is a subgraph whose vertex set is the entire vertex set of $G$.
We refer the reader to \cite{Bondy} for the terminology and notations not defined in this paper.

By the definition of the tree edge-connectivity, the following result is obvious.

\begin{obv}[\cite{LWMM}]\label{obv1}
For any integer $2\leq s\leq t$, $\lambda_s(G)\geq \lambda_t(G)$.
\end{obv}

Now, we give a confirmative solution to Conjecture \ref{conj1}.

\begin{thm}\label{thm1}
Let $k\geq 2$ be an integer. If a connected graph $G$ has at least $k$ vertices and at least $k$ edges, then
$\kappa_k(L(G))\geq \lambda_k(G)$. Moreover, the bound is sharp.
\end{thm}
\pf Let $v_e$ be the vertex of the line graph $L(G)$ corresponding to the edge $e$ of $G$.
Assume that $\lambda_k(G)=m$. For any $k$-subset $S_L=\{v_{e_1},v_{e_2},\ldots,v_{e_k}\}$ of $V(L(G))$,
by the definition of $\kappa_k$, it suffices to show that $\kappa_{L(G)}(S_L)\geq m$.

Now, let $S_G=\{e_1,e_2,\ldots,e_k\}$ and then $S_G\subseteq E(G)$. Denote by $G[S_G]$
the edge-induced subgraph of $G$ whose edge set is $S_G$ and whose
vertex set consists of all ends of edges of $S_G$.

We distinguish two cases:

{\bf Case 1:} None of the components of $G[S_G]$ is a tree or unicyclic.

Therefore, for each component $C_l$ of $G[S_G]$, $|E(C_l)|-|V(C_l)|\geq 1$ and so
$|E(G[S_G])|=|S_G|=k\geq |V(G[S_G])|+1$. Let $V(G[S_G])=Q^*$. It follows that
$|Q^*|\leq k-1$. Since $G$ has at least $k$ vertices, we can take a vertex $v^*$ in
$V(G)\setminus Q^*$ and then let $Q=Q^*\cup \{v^*\}$. Now, because $|Q|\leq k$,
by Observation \ref{obv1}, there are $m$ edge-disjoint
$Q$-Steiner trees $T_1,T_2,\ldots,T_m$ in $G$.

Next, in each tree $T_r$ ($1\leq r\leq m$), we will assign a specific edge to each
vertex of $Q^*$. To see this, we let $v^*$ be the root and define the level $l(v)$
of a vertex $v$ in $T_r$ to be the distance from the root $v^*$ to $v$.
It is easy to see that, for each vertex $v_i\in Q^*=Q\setminus \{v^*\}$,
there is a unique edge $e$ connecting the vertex $v_i$ and a vertex of level $l(v_i)-1$.
Assign the edge $e$ to the vertex $v_i$. We say that the edge $e$ is the corresponding edge
of $v_i$ in $T_r$ and denoted by $\hat{e}^r_i$. Note that any two vertices of $Q^*$ in $T_r$
have different corresponding edges. More precisely, $\hat{e}^r_i\neq \hat{e}^r_j$
for any $1\leq i\neq j\leq |Q^*|$.

Now, for each tree $T_r$ ($1\leq r\leq m$) and each edge $e=v_iv_j\in S_G$,
do the following operation. Note that, by Lemma \ref{lem1},
$L(T_r)$ ($1\leq r\leq m$) is a connected subgraph of $L(G)$. Moreover,
since $Q^*=V(G[S_G])$, both ends of each edge in
$S_G$ belong to $Q^*$ and so $v_i,v_j\in V(T_r)$.

\noindent {\bf Operation A:} If $e\in E(T_r)$, it is done; otherwise $e\notin E(T_r)$,
that is $v_e\notin V(L(T_r))$. Note that, $T_1,T_2,\ldots,T_m$ are edge-disjoint
and so at most one of them contains the edge $e$.

If $e\in E(T_s)$, where $1\leq s\neq r\leq m$, then $e$ is the corresponding edge
of one of its ends in $T_s$. Without loss of generality, assume that $e$ is the
corresponding edge of $v_i$ in $T_s$, that is, $\hat{e}^s_i=e$.
Now, for $T_r$, there is an edge $\hat{e}^r_j$ corresponding to the vertex $v_j$,
which is the other end of $e$. Since $e$ and $\hat{e}^r_j$ have the same end $v_j$,
they are adjacent and so $v_{e}v_{\hat{e}^r_j}\in E(L(G))$.
Add the vertex $v_{e}$ and the edge $v_{e}v_{\hat{e}^r_j}$ to $L(T_r)$.

Otherwise, none of the trees $T_1,T_2,\ldots,T_m$ contains the edge $e=v_iv_j$.
In this case, we can add the vertex $v_{e}$ and either the edge $v_{e}v_{\hat{e}^r_i}$
or the edge $v_{e}v_{\hat{e}^r_j}$ to $L(T_r)$, where $\hat{e}^r_i$ and $\hat{e}^r_j$
are the corresponding edges of $v_i$ and $v_j$ in $T_r$, respectively.     $\square$

Now, $L(T_1),L(T_2),\ldots,L(T_m)$ are transformed into $m$ connected subgraphs of $L(G)$,
each of which contains the vertex set $S_L$. Next, for each of the obtained subgraphs of $L(G)$, take a
spanning tree $T^*_r$ ($1\leq r\leq m$). Because $V(T^*_r)\supseteq S_L$ ($1\leq r\leq m$), it remains
to show that $T^*_1,T^*_2,\ldots,T^*_m$ are internally disjoint. Note that, if $v_{e}\notin V(L(T_r))$,
for some $e\in S_G$ and $r\in \{1,2,\ldots,m\}$, $v_{e}$ must be a leaf of $T^*_r$.

Since $T_1,T_2,\ldots,T_m$ are edge-disjoint in $G$, $L(T_1),L(T_2),\ldots,L(T_m)$ are
vertex-disjoint in $L(G)$. Moreover, the vertices added to $L(T_r)$ by Operation A are all from $S_L$.
Therefore, $T^*_1,T^*_2,\ldots,T^*_m$ are vertex-disjoint except $S_L$, that is, $V(T^*_r)\cap V(T^*_s)=S_L$,
for any $1\leq r<s\leq m$.

Now, assume that there are two trees $T^*_r$ and $T^*_s$ such that $E(T^*_r)\cap E(T^*_s)\neq \emptyset$
($1\leq r< s\leq m$). Let $e^* \in E(T^*_r)\cap E(T^*_s)$. Since  $V(T^*_r)\cap V(T^*_s)=S_L$,
both ends of $e^*$ belong to $S_L$. Thus, without loss of generality, let $e^*=v_{e_1}v_{e_2}$. If $L(T_r)$
contains neither $v_{e_1}$ nor $v_{e_2}$, by Operation A, both $v_{e_1}$ and $v_{e_2}$ are leaves of $T^*_r$
and hence it is impossible that $v_{e_1}v_{e_2}\in E(T^*_r)$. So is $L(T_s)$. And
$L(T_r)$ and $L(T_s)$ are vertex-disjoint ($T_r$ and $T_s$ are edge-disjoint). Thus,
without loss of generality, suppose that $v_{e_2}\in L(T_r)$ and $v_{e_1}\in L(T_s)$,
and so $v_{e_1}\notin L(T_r)$ and $v_{e_2}\notin L(T_s)$.
Since $v_{e_1}$ and $v_{e_2}$ are adjacent in $L(G)$, $e_1$ and $e_2$ are adjacent in $G$.
Assume that $v_i$ is the common end of $e_1$ and $e_2$ in $G$ and let $e_1=v_iv_j$.
Since $v_{e_2}v_{e_1} \in E(T^*_s)$, we added the vertex $v_{e_2}$ and the edge $v_{e_2}v_{e_1}$
to $L(T_s)$. So by Operation A, we know that $e_1$ is exactly the corresponding edge of $v_i$
in $T_s$, that is, $e_1=\hat{e}^s_i$. Again by Operation A, since $e_1\notin E(T_r)$,
we added the vertex $v_{e_1}$ and the edge $v_{e_1}v_{\hat{e}^r_j}$ to $L(T_r)$, where
$e_1$ and $\hat{e}^r_j$ have the same end $v_j$. Since $e_1\neq e_2$ and
$e_1$ and $e_2$ have the same end $v_i$, it is impossible that $\hat{e}^r_j=e_2$.
Therefore, by Operation A, it is impossible that $v_{e_1}v_{e_2}=e^* \in E(T^*_r)$,
a contradiction. It follows that $T^*_1,T^*_2,\ldots,T^*_m$ are edge-disjoint.

Thus, in this case, $T^*_1,T^*_2,\ldots,T^*_m$ are $m$ internally disjoint trees connecting $S_L$ in $L(G)$.

{\bf Case 2:} There is a component of $G[S_G]$ which is either a tree or unicyclic.

For each component $C_l$ of $G[S_G]$ which is neither a tree nor unicyclic, add
all vertices of $C_l$ to the vertex set $Q_1$ and add all edges of $C_l$ to the edge set $S^1_G$.
Clearly, if $Q_1\neq \emptyset$, $|S^1_G|> |Q_1|$.

Next, for each component $C_t$ of $G[S_G]$ which is either a tree or unicyclic, if
$C_t$ is unicyclic, choose an edge $e_t$ from $C_t$ such that $C_t-e_t$ is a tree
and let one end of $e_t$ as the root $r_t$; otherwise, select an arbitrary vertex
as the root $r_t$. For $C_t$ (if $C_t$ is a tree) or $C_t-e_t$ (if $C_t$ is unicyclic),
define the level $l(v)$ of a vertex $v$ to be the distance from the root $r_t$ to $v$.
Notice that each edge in the tree $C_t$ (or $C_t-e_t$ if $C_t$ is unicyclic) joins
vertices on consecutive levels. Then, for each edge $e=uv$, where $l(u)+1=l(v)$,
we assign the vertex $v$ which has higher level to the edge $e$ and say that the
vertex $v$ is the corresponding vertex of the edge $e$. If $C_t$ is unicyclic,
let the root $r_t$ be the corresponding vertex of the remaining edge $e_t$.
Now, each edge of $C_t$ has a corresponding vertex. By the definition,
it is obvious that any two edges of $C_t$ have different corresponding vertices.
Add the corresponding vertices of all edges of $C_t$ to the vertex set $Q_2$
and add all edges of $C_t$ to the edge set $S^2_G$. Clearly, $|S^2_G|=|Q_2|$.

Moreover, it is clear that $Q_1\cap Q_2=\emptyset$, $S^1_G\cap S^2_G=\emptyset$
and $S_G=S^1_G\cup S^2_G$. Let $S^1_L=\{v_{e}| e\in S^1_G\}$ and
$S^2_L=\{v_{e}| e\in S^2_G\}$. Then $S_L=S^1_L\cup S^2_L$.
Let $Q=Q_1\cup Q_2$. We have $|Q|=|Q_1|+|Q_2|\leq |S^1_G|+|S^2_G|=|S_G|=k$.
Since $Q\subseteq V(G)$, by Observation \ref{obv1}, there are $m$ edge-disjoint
$Q$-Steiner trees $T_1,T_2,\ldots,T_m$ in $G$. Note that both ends of each edge in $S^1_G$
belong to $Q_1$, but there may be an edge in $S^2_G$, only one end of which
belongs to $Q_2$. Thus, we use different methods to deal with the edges in $S^1_G$
and $S^2_G$.

For every edge of $S^1_G$, we take the same approach as Case 1.
In each tree $T_r$ ($1\leq r\leq m$), since $Q_2\neq \emptyset$ (it is possible
that $Q_1=\emptyset$), take an arbitrary vertex $v^*$ in $Q_2$ as the root
and define the level $l(v)$ of a vertex $v$ in $T_r$ to be the distance from the root
$v^*$ to $v$. For each vertex $v_i\in Q_1$ (if $Q_1\neq \emptyset$),
there is a unique edge $e$ connecting the vertex $v_i$ and a vertex of level $l(v_i)-1$.
Let the edge $e$ be the corresponding edge of $v_i$ in $T_r$, denoted by $\hat{e}^r_i$.
Any two vertices of $Q_1$ in $T_r$ have different corresponding edges.
Now, apply Operation A to each tree $T_r$ ($1\leq r\leq m$) and each edge $e=v_iv_j\in S^1_G$.
Then, each vertex of $S^1_L$ is added to $L(T_r)$ ($1\leq r\leq m$).

Next, for each edge $e_i$ of $S^2_G$ and each tree $T_r$ ($1\leq r\leq m$),
do the following operation.

\noindent {\bf Operation B:} If $e_i\in E(T_r)$, it is done; otherwise $e_i\notin E(T_r)$,
that is $v_{e_i}\notin V(L(T_r))$. By the definitions of $S^2_G$ and $Q_2$,
there is a corresponding vertex $v_i$ of $e_i$, and $v_i\in Q_2\subseteq Q$ and so $v_i\in V(T_r)$.
Thus, there exists an edge $\tilde{e}^r_i\neq e_i$ incident with $v_i$ in the tree $T_r$.
Since $e_i$ and $\tilde{e}^r_i$ have the same end $v_i$, they are adjacent
and so $v_{e_i}v_{\tilde{e}^r_i}\in E(L(G))$.
Add the vertex $v_{e_i}$ and the edge $v_{e_i}v_{\tilde{e}^r_i}$ to $L(T_r)$.

Now, after applying Operations A and B, $L(T_1),L(T_2),\ldots,L(T_m)$ are
transformed into $m$ connected subgraphs of $L(G)$,
each of which contains the vertex set $S_L=S^1_L\cup S^2_L$.
For each of the obtained subgraphs of $L(G)$, take a spanning tree $T^*_r$ ($1\leq r\leq m$).
Note that, if $v_{e}\notin V(L(T_r))$, for some $e\in S_G$ and $r\in \{1,2,\ldots,m\}$,
whether $e\in S^1_G$ or $S^2_G$, that is, whether Operation A or Operation B is applied,
$v_{e}$ must be a leaf of $T^*_r$.
Because $V(T^*_r)\supseteq S_L$ for any $1\leq r\leq m$, it remains
to show that $T^*_1,T^*_2,\ldots,T^*_m$ are internally disjoint.

Since $L(T_1),L(T_2),\ldots,L(T_m)$ are vertex-disjoint in $L(G)$ and
the vertices added to $L(T_r)$ by Operations A and B are all from $S_L$,
$T^*_1,T^*_2,\ldots,T^*_m$ are vertex-disjoint except $S_L$, that is, $V(T^*_r)\cap V(T^*_s)=S_L$,
for any $1\leq r<s\leq m$.

To complete the proof, it remains to show that $T^*_1,T^*_2,\ldots,T^*_m$ are edge-disjoint.
By contradiction, assume that there are two trees $T^*_r$ and $T^*_s$
such that $E(T^*_r)\cap E(T^*_s)\neq \emptyset$
($1\leq r< s\leq m$). Let $e^* \in E(T^*_r)\cap E(T^*_s)$. Since  $V(T^*_r)\cap V(T^*_s)=S_L$,
both ends of $e^*$ belong to $S_L$. Thus, without loss of generality, let $e^*=v_{e_1}v_{e_2}$. If $L(T_r)$
contains neither $v_{e_1}$ nor $v_{e_2}$, then whether apply Operation A or Operation B,
both $v_{e_1}$ and $v_{e_2}$ are leaves of $T^*_r$, which is impossible.
So is $L(T_s)$. And $L(T_r)$ and $L(T_s)$ are vertex-disjoint ($T_r$ and $T_s$ are edge-disjoint). Thus,
without loss of generality, suppose that $v_{e_2}\in L(T_r)$ and $v_{e_1}\in L(T_s)$,
and so $v_{e_1}\notin L(T_r)$ and $v_{e_2}\notin L(T_s)$.
Since $v_{e_1}$ and $v_{e_2}$ are adjacent in $L(G)$, $e_1$ and $e_2$ are adjacent in $G$.
Therefore, $e_1$ and $e_2$ belong to the same component of $G[S_G]$.
Hence, by the definitions of $S^1_G$ and $S^2_G$, both $e_1$ and $e_2$
belong to $S^1_G$ or $S^2_G$.

If both $e_1$ and $e_2$ belong to $S^1_G$, by Operation A, it is
impossible that $v_{e_1}v_{e_2}=e^*\in E(T^*_r)\cap E(T^*_s)$. The proof is
the same as that of Case 1.

If both $e_1$ and $e_2$ belong to $S^2_G$, since $e_1\notin E(T_r)$,
by Operation B, we added the vertex $v_{e_1}$ and the edge $v_{e_1}v_{\tilde{e}^r_1}$ to $L(T_r)$,
where the common end $v_1$ of $e_1$ and $\tilde{e}^r_1$ in $G$ is the
corresponding vertex of $e_1$. Similarly, since $e_2\notin E(T_s)$,
by Operation B, we added the vertex $v_{e_2}$ and the edge $v_{e_2}v_{\tilde{e}^s_2}$ to $L(T_s)$,
where the common end $v_2$ of $e_2$ and $\tilde{e}^s_2$ in $G$ is the
corresponding vertex of $e_2$. Since $v_1\neq v_2$ by the definition of $Q_2$,
at least one of the equations $\tilde{e}^r_1=e_2$ and $\tilde{e}^s_2=e_1$
is not true. So $v_{e_1}v_{\tilde{e}^r_1}\neq v_{e_1}v_{e_2}$ or
$v_{e_2}v_{\tilde{e}^s_2}\neq v_{e_1}v_{e_2}$.
It is impossible that $e^*=v_{e_1}v_{e_2}\in E(T^*_r)\cap E(T^*_s)$, a contradiction.
It follows that $T^*_1,T^*_2,\ldots,T^*_m$ are edge-disjoint.

Thus, in both cases, there always exist $m$ internally disjoint trees connecting $S_L$ in $L(G)$
and so $\kappa_{L(G)}(S_L)\geq m$. By the arbitrariness of $S_L$, we
conclude that $\kappa_k(L(G))\geq m$.

For a cycle $C_n$ with $n\geq k$, since $L(C_n)=C_n$,
$\kappa_k(L(C_n))=\lambda_k(C_n)=1$ for $k\geq 3$ and $\kappa_2(L(C_n))=\lambda_2(C_n)=2$.
Thus, the bound is sharp.
The proof is complete.\qed

\section*{Acknowledgments.}

The author's work was supported by Zhejiang Provincial Natural Science Foundation of China (No. LY18A010002)
and the Natural Science Foundation of Ningbo, China.

\end{document}